\documentclass[francais]{smfart} 
\usepackage{correspondance-c-w}

\begin{document}
\title[Une lettre d'Henri Lebesgue à \'Elie Cartan]{Une lettre d'Henri Lebesgue à \'Elie Cartan}
\author{Michèle Audin}
\address{Institut de Recherche Mathématique Avancée\\
Université de Strasbourg et CNRS\\
7 rue René Descartes\\
67084 Strasbourg\index{Strasbourg} cedex\\
France}
\email{Michele.Audin@math.u-strasbg.fr}
\urladdr{http://www-irma.u-strasbg.fr/~maudin}
\thanks{Version du \today}

\begin{abstract}
Nous publions une lettre de Lebesgue à \'Elie Cartan, ainsi qu'une lettre de Montel à \'Elie Cartan, datant de la fin 1933, où il est question de Gaston Julia, de Paul Montel et d'une élection à l'Académie des sciences. Nous discutons le contexte historique et les mathématiques.
\end{abstract}

\begin{altabstract}
We publish a letter from Lebesgue to Cartan and a letter from Montel to Cartan, dated 1933--1934, about Gaston Julia, Paul Montel, and an election at the Paris Academy of Sciences. We discuss the context and the mathematics.
\end{altabstract}



\maketitle

\subsection*{Avertissement}

Les deux lettres ont été découvertes entre la publication de~\cite{GPSM18} et de sa version en anglais~\cite{GPSM18eng}. Elles ont donc été incluses dans la traduction. Le but de cet article est d'en faire profiter les lecteurs de la version française.

\Section*{Introduction}
C'est en 1988 qu'a été découvert, dans les sous-sols de l'Institut Henri Poincaré, le trésor des lettres conservées par son premier directeur \'Emile Borel. Parmi celles-ci, de nombreuses lettres d'Henri Lebesgue, qui ont été publiées dans le \emph{Cahier du séminaire d'histoire des mathématiques}~\cite{LebesgueBorel}, avant de faire l'objet du livre~\cite{LebesgueBorel2}. Dans~\cite{LebesgueBorel}, un \og article\fg tapé à la machine, sans aucun système de numérotation automatique, dont les notes ont été rédigées et ajoutées par Pierre Dugac au cours de la frappe, on peut lire successivement
\begin{itemize}
\item Note [900], 
\begin{quote}
\noindent Henri Cartan nous a dit le 1er octobre 1990 qu'il possédait des lettres de Lebesgue à son père.
\end{quote}
\item Note [977], 
\begin{quote}
\noindent Henri Cartan a bien voulu nous envoyer la seule lettre qu'il possède de Lebesgue à son père, datée du 14 avril 1930 (cachet de la poste).
\end{quote} 
Dans cette lettre, qui est reproduite dans la note en question, Lebesgue demande à Cartan d'examiner la thèse de Georges de Rham.
\item Note [1064],
\begin{quote}
\noindent Henri Cartan a bien voulu nous remettre encore une lettre de Lebesgue à son père qu'il venait de retrouver.
\end{quote} 
Dans cette lettre, datée du 19 octobre 1925 et reproduite intégralement elle aussi par Dugac, Lebesgue recommande Raoul Dautry à Cartan, puis il répond à une question que Cartan lui avait posée deux ans auparavant, sur le plan projectif complexe \og du point de vue de l'analysis situs\fg, il en donne les nombres de Betti (sous la forme $BP_1=BP_3=1$, $BP_2=2$ (\string?) et s'étonne que le groupe fondamental de cette variété soit trivial.
\end{itemize}
Lors de la publication du livre~\cite{LebesgueBorel2}, la situation avait encore évolué, puisque l'on peut lire dans ce livre:
\begin{itemize}
\item Note [555],
\begin{quote}
\noindent Henri Cartan nous a indiqué qu'il possède une lettre de Lebesgue à \'Elie Cartan où Lebesgue l'attrape parce qu'il a soutenu Julia contre Montel.
\end{quote}
\end{itemize}

C'est cette dernière lettre qui fait l'objet du présent article. Notre intérêt pour elle et son contenu vient des recherches ayant mené à l'écriture du livre~\cite{GPSM18}: les relations entre Gaston Julia et Paul Montel, et notamment le moment, crucial pour ces relations, de l'élection de Julia à l'Académie des sciences en 1934, jouent en effet un rôle important dans ce livre.

La lettre en question se trouvait, Henri Cartan l'a dit à Dugac, à son domicile, 95 boulevard Jourdan. Cet appartement spacieux équipé de nombreux placards et armoires, a été occupé par \'Elie Cartan (et sa famille de) 1937 à sa mort en 1951, puis par Henri Cartan et sa famille jusqu'en 2008. Nous avons retrouvé cette lettre, aujourd'hui aux archives de l'Académie des sciences, en aidant la famille Cartan à ranger les gigantesques archives mathématiques d'\'Elie et Henri Cartan.

Avant d'en livrer le texte et de laisser la parole à Lebesgue, avec son style personnel, vivant et... inimitable, rappelons-en le, les contextes. 

\section{Contextes}
Le contexte immédiat est celui de l'élection, à l'Académie des sciences et dans la section de géométrie, d'un nouveau membre, en remplacement de Paul Painlevé, décédé le 29 octobre 1933. Cette élection aura lieu le 5 mars 1934. C'est Gaston Julia qui sera élu.

Un contexte un peu plus général est celui des relations scientifiques et de pouvoir entre les mathématiciens français dans cette période qui n'était pas encore l'Entre-deux-guerres, mais était, de façon inoubliable, une Après-guerre. Les deux personnages principaux que la lettre de Lebesgue met en scène sont Gaston Julia (1893--1978) et Paul Montel (1876--1975). 

De leur histoire scientifique, rappelons que
\begin{itemize}
\item Paul Montel est l'inventeur des familles normales de fonctions holomorphes, et en particulier a eu l'idée d'appliquer ses résultats sur les \emph{suites} de fonctions à la démonstration par exemple de résultats sur \emph{une} fonction (les théorèmes de Picard notamment).
\item Gaston Julia, plus jeune, est un grand blessé de la grande guerre (une \og gueule cassée\fg) et l'un des rares survivants mathématiciens de sa génération. Il a obtenu deux prix de l'Académie des sciences pendant la guerre, le prix Bordin pour sa thèse en 1917 et le Grand prix des sciences mathématiques en 1918 pour ses travaux sur l'itération.
\end{itemize}

Il est connu, et a d'ailleurs été reconnu par les deux auteurs eux-mêmes, que les travaux sur l'itération de Julia et de Fatou ont été stimulés par la parution d'une note de Montel qui leur a démontré l'usage qu'ils pouvaient faire des familles normales de fonctions (sur cette histoire, voir~\cite{GPSM18} et le \T\ref{subsecmathmontel} ci-dessous). Rappelons aussi que Fatou, qui avait renoncé à concourir pour le Grand prix de 1918, est mort en 1929. Nous ne le verrons donc pas apparaître ici.

\medskip 
Il est connu que Lebesgue a très mal vécu les arguments utilisés lors de cette élection, au point qu'il a arrêté d'assister aux séances de l'Académie des sciences. Un passage un peu cryptique des souvenirs de Lucienne Félix~\cite{FelixLebesgue}
\begin{quote}
\`A l'Académie [...] il s'indigna qu'on tint compte pour les élections, d'autres éléments que la valeur scientifique [...] Il alla, paraît-il, jusqu'à tenir tête au grand Maître et Secrétaire perpétuel \'Emile Picard, et pendant un certain temps, cessa d'aller aux Séances.\end{quote}
éclairé par Montel lui-même~\cite{Montelanecdotes}
\begin{quote}
Henri Lebesgue\index{Lebesgue (Henri), 1875--1941, mathématicien}, grand mathématicien du début de ce siècle; pour moi, l'ami le plus généreux, le plus sûr et le plus loyal.

En 1934, je posais ma candidature à l'Académie des Sciences. Je ne fus pas élu, les suffrages étant allés à un grand mutilé de la guerre de 14. Lebesgue considéra qu'une injustice avait été commise et refusa, pendant trois ans, de paraître aux séances.

Il n'y retourna qu'en 1937, après mon élection, qui fut obtenue par un nombre élevé de voix [...]
\end{quote}
nous le rappelle. Cette lettre nous montre Lebesgue avant l'élection, avec sa personnalité entière, mais surtout avec ses arguments scientifiques et déontologiques.

\medskip
Dans la suite de cet article, nous discutons brièvement la question de la date de la lettre de Lebesgue à Cartan, nous présentons le texte de cette lettre (avec de brèves notes explicatives) ainsi que celui d'une lettre de Montel à Cartan datant de cette période (et datée, elle), puis nous reprenons la parole pour des commentaires plus précis et pour  évoquer l'élection et sa suite.

\section{Lettres de Lebesgue et Montel à \'Elie Cartan}
Les deux lettres sont manusctites. Celle de Lebesgue est assez longue (quatre feuilles numérotées de 1 à 4). Nous indiquons par des / les changements de page.

\subsection*{Les dates de ces lettres}
La lettre de Lebesgue n'est pas datée. Et ce n'est pas le \og Lundi à ma rentrée de l'Institut\fg qui peut aider à en préciser la date: les séances ordinaires de l'Académie des sciences ont lieu le lundi! Il y est question d'une élection à l'Académie des sciences où s'affrontent Julia et Montel et il n'y en a eu qu'une, puisque Julia a été élu dès son coup d'essai (qui fut un coup de maître). Novembre, décembre 1933 ou janvier 1934, donc. Plus vraisemblablement décembre ou janvier, puisque Lebesgue mentionne \og le début de novembre\fg (et pas \og le début du mois\fg). 

La lettre de Montel que nous publions aussi est datée, elle, du 24 décembre 1933. Elle se situe exactement dans la même période, au sein de la même discussion. Sa position par rapport à la lettre de Lebesgue n'est pas absolument claire: soit les éclaircissements donnés par Montel font partie de l'\og épluchage\fg mentionné par Lebesgue et la lettre est postérieure, soit la discussion avec Montel et sa lettre ont été provoquées par la lettre de Lebesgue et elle est légèrement antérieure...

\subsection*{La lettre de Lebesgue}\mbox{ }

\begin{flushright}
Mercredi matin
\end{flushright}

\medskip
Mon cher Cartan

\medskip
C'est du lit que je vous réponds; je me suis couché Lundi à ma rentrée de l'Institut brûlant de fièvre, rechute de cette grippe que je traîne depuis le début de novembre sans avoir pu m'en libérer. Le mardi ma femme est venue me retrouver et nous voici cote à cote toussant, éternuant, crachant, ce qui est infiniment touchant.

Ceci m'a empêché de convoquer Julia. Votre bienveillance à vous est infinie et après qu'il vient de vous faire passer 15 jours à éplucher mot par mot les travaux de Montel pour prouver que ses travaux à lui ne doivent rien, ou si peu, à ceux de Montel --- à éplucher comme jamais vous n'avez épluché aucun écrit, à éplucher d'une façon telle que si l'on faisait le même travail sur Poincaré il n'y aurait rien, absolument rien de démontré par Poincaré --- vous vous satisfaites, précisément dans ce mémoire qu'il s'agit de montrer indépendant de Montel, d'une citation où il est dit, je crois, que Montel a généralisé le Lindelöf et quand après cela le théorème cité est celui de Lindelöf inopérant et non celui de Montel utilisé. Après cela aussi Montel devient de l'Ecole scandinave et la première méthode, étant scandinave, est indépendante de Montel!\string!\string!\string!

Et pourtant, vous m'avez concédé que Julia a \og\emph{amenuisé}\footnote{Selon une tradition établie, les italiques rendent ici les soulignements du texte manuscrit. Toutes les notes sont de l'auteur.}\fg le rôle de Montel dans sa fameuse communication de Zurich\footnote{Il s'agit de la conférence plénière prononcée par Julia pensant le congrès international de Zurich en 1932. Julia était très fier de ce texte~\cite{Julia32}, qu'il a fait publier par Gauthier-Villars en volume séparé.}!

Non, en vérité! Cartan, vous exagérez en ne voulant pas voir ce qui est éclatant. Que vous en excusiez Julia, cela c'est tout autre chose; il est à part et même si on reconnaît en lui un travers ou même une tare on peut très légitimement prétendre que c'est une déformation mentale due à l'état lamentable dans lequel la guerre l'a laissé et que vous ne l'en plaigniez que davantage je n'aurais rien à y objecter /
car c'est à cette mansuétude que je voudrais arriver. Vous ne m'y aidez pas; en contestant l'évidence vous me faites [? ranimer, raviver] davantage les reproches qu'on lui peut faire. Et j'aurais pourtant bien besoin d'être aidé car rien n'est plus difficile que de pardonner à ceux qui sont tellement personnels qu'ils sont injustes envers les autres car rien n'est plus contraire à mon tempérament. Et je n'ai jamais vu personne d'aussi âprement personnel que Julia.

Permettez-moi de vous rappeler ma façon d'agir, elle est tellement différente de celle de Julia qu'elle vous expliquera la difficulté que j'ai à encaisser ses façons de faire. Quand certains ont cherché à me pousser dans les jambes de Borel, je me suis dit, qu'aurais-tu fait si Borel n'avait pas existé? Ma vanité me dictait pas mal de réponses \og avantageuses\fg je ne l'ai pas écoutée et je me suis dit: tu ne peux répondre avec certitude, tu dois laisser passer Borel mais même exprimer, ce qui est la vérité, c'est qu'il est logiquement, tout aussi bien que chronologiquement, avant toi. Et je l'ai fait sans aucune réticence. J'y avais peut-être quelque mérite, car Borel ne m'avait pas toujours traité avec justice; mais ce que j'avais eu à lui reprocher, je n'avais pas été le dire à l'oreille, je le lui ai dit crûment et publiquement, sans aucune diplomatie, sans ménagement, mais aussi en toute sincérité. Je défie qu'on trouve dans mon âcre revendication un mot injuste pour l'\oe uvre de Borel, un essai tendancieux de me pousser à son détriment. Et ma philippique avait assez porté pour qu'elle m'ait puissamment aidé dans une campagne contre Borel si j'avais voulu en faire une; j'avais pour moi Goursat et dans la section, Picard Koenigs en dehors d'elle\footnote{Il s'agit de la Section de géométrie qui regroupait (la suite de cette phrase est un anachronisme) les mathématiciens \og purs\fg de l'Académie des sciences. Koenigs était membre de la Section de mécanique. Rappelons qu'il y avait des mathématiciens dans les sections de géométrie, de mécanique et d'astronomie.} et qui ne demandaient qu'à faire campagne, j'ai dit non. Et à moi je me suis dit: une place /

(2)

\noindent à l'Académie ne vaut pas qu'on se diminue à ses propres yeux, tu attendras.

Et quand Jordan est mort je n'ai fait aucune démarche; je n'ai été voir personne, ni dans la Section, ni en dehors d'elle. C'est seulement quand la Section m'a dit: nous vous présentons à l'unanimité, que j'ai commencé à rédiger ma Notice; si bien qu'il a fallu reculer l'élection. Cette élection a été ce que vous savez; tout le monde y a eu des voix, sauf vous\footnote{Jordan est mort le 21 janvier 1922 et l'élection a eu lieu un peu plus de quatre mois plus tard, le 29 mai. Et il n'est pas tout à fait vrai que Cartan n'y a eu aucune voix. Voir le~\T\ref{secelection}.\label{notemortjordan}}. Et cela ne m'a pas étonné car j'avais pu constater au cours des visites que deux personnes seulement ne demandaient pas à ceux qu'ils visitaient de leur donner leur voix: moi, qui n'y avais guère de mérite, et vous.

J'avais [aurais?] cru que, dans la circonstance actuelle, vous condamneriez comme moi tout acte révélant une ambition poussée au point de ne pas s'arrêter devant l'injustice puisque vous aviez su attendre, plus encore que moi, que le jugement des autres vous soit favorable, quelque jugement personnel que vous puissiez légitimement avoir. Et je répète que condamner les procédés, ce n'est pas nécessairement condamner l'homme quand il s'agit d'un grand mutilé.

Comme nous voici loin de la sérénité et du genre de préoccupations purement scientifiques qui seules devraient nous importer dans une élection scientifique. De qui est-ce la faute? De Montel?

De quoi s'agit-il? Denjoy s'étant retiré devant Montel, ce qui tout de même est un hommage de poids envers celui-ci, il faut savoir si Julia et Montel sont dignes de l'Académie et dans quel ordre ils doivent y être appelés. Si nous jugeons qu'ils en sont tous deux dignes, quel que soit l'ordre dans lequel nous les placions scientifiquement --- et je me hâte / de répéter que je place Montel en premier --- une différence de 17 ans d'âge compte et est péremptoire humainement mais ce qui pour moi compte bien autrement c'est que nul ne pourrait dire ce que Julia aurait fait si Montel n'avait pas existé car ses deux travaux les plus importants relèvent de Montel.

Je suis scientifiquement pour Montel parce que là où tout le monde marchait au hasard et sans bien comprendre ce qu'il faisait il a fait voir à tous le fait fondamental duquel tout découlait et sur lequel devait porter l'effort. Des assortiments de théorèmes sont devenus des cas particuliers d'un même fait bien simple. Là où la veille il fallait de l'invention, de l'ingéniosité et du bonheur, il n'a plus fallu le lendemain que l'application d'une méthode. Cela, c'est une très grande chose. Puisque l'un des mérites de cette chose est de faire comprendre ce qui l'a précédée, il faut qu'il y ait des prédécesseurs. On a pu prétendre que la géométrie analytique de Descartes et Fermat ne faisait que répéter ce qui était d'usage courant chez les mathématiciens depuis Apollonius; que la dérivation de Newton était de Barrow, de Fermat, etc. Tout cela ne me trouble pas. 

Aussi parce qu'à côté de cette preuve d'intelligence qui permet de voir les faits de haut il a tout de même fait bien d'autres choses --- on a dans l'examen des pointes d'aiguilles soulevées par Julia trop tendance à l'oublier --- et par les familles normales et autrement. En particulier je dirai que j'apprécie fort les difficultés vaincues par Montel pour arriver, enfin, à montrer que les conditions de Cauchy 
$$\dfrac{\partial P}{\partial x}=\dfrac{\partial Q}{\partial y}, \dfrac{\partial P}{\partial y}=-\dfrac{\partial Q}{\partial x}$$ 
sont suffisantes, sans rien de plus. Ses recherches sur l'existence des dérivées des fonctions d'une ou de plusieurs variables sont aussi de ma compétence et ce dernier cas n'était pas bien abordé avant lui.

Je ne veux pas faire une énumération, inutile avec vous qui êtes la conscience même; il n'était sans doute même pas nécessaire de vous rappeler que Montel n'est pas seulement les familles normales.

Julia a à son actif quantité de beaux travaux mais je crois /

(3)

\noindent
qu'il restera moins de ses travaux que de ceux de Montel, qu'ils seront moins utiles aux progrès de la science que ceux de Montel ne l'ont déjà été. Pourquoi? parce que, hors les droites de Julia qui sont de beaucoup son principal titre, il n'a jamais été un initiateur, que dans ses travaux, le foisonnement, l'érudition, la technique savante, parfois inutilement, nous impressionnent et que tout cela est en réalité extérieur à la vraie valeur. Cette valeur est réelle et je ne pense nullement ni à la contester, ni à l'\og amenuiser\fg mais il fallait que je vous dise ce qui influe sur mon jugement.

Reste à savoir si, oui ou non, les deux principaux travaux de Julia utilisent comme outils essentiels les outils créés par Montel?  Pour l'itération, Julia, Fatou et le rapporteur Humbert l'ont reconnu formellement mais cet aveu pèse peut-être à Julia. Moi, qui manque de votre bienveillance infinie, j'imagine qu'il vous a dit: j'étais mal renseigné, j'aurais pu et dû citer Landau et Carathéodory et non Montel --- car, tout de même vous n'avez pas trouvé seul Landau et Carathéodory. C'est bien assez de lire Montel et Julia s'il fallait encore lire tous ceux qu'ils citent, où irions nous?

Donc, alerté, vous avez lu Landau et Carathéodory et y avez trouvé ce que vous m'avez dit. Imaginons que Landau et Carathéodory aient non seulement construit incidemment une démonstration voisine de celle du critère de Montel ou cette démonstration même mais qu'ils en aient eu conscience et qu'ils l'aient énoncé; eh bien? \c Ca empêchera-t-il que ce qu'utilise Julia ne provienne des idées de Montel puisque, m'avez-vous dit, Landau et Carathéodory se réfèrent expressément aux idées et résultats de Montel. Donc sur ce premier point (et quoiqu'il arrive de la question \emph{différente} priorité Montel-Landau) j'ai une opinion ferme.

(Quant à la question de priorité, demandez à Montel. Je ne suis pas très impressionné par la question; un mémoire paru dans un périodique \emph{daté} 1911 et un résultat annoncé dans une note de 1911 ne peuvent pas avoir beaucoup d'influence l'un sur l'autre. Le mémoire se réfère à Montel voila le critère et Montel connaissait le mémoire lors de son article développé de 1912; mais ce peuvent être deux connaissances assez différentes) /

Pour les droites de Julia. Triple offensive de Julia:
\begin{enumerate}
\item[a] Points $J$. Liquidée n'est-ce pas? voir le 1\up{er} paragraphe du mémoire Montel 1912, reproduit presque textuellement (avec bien d'autres choses du même mémoire) dans le livre de Julia. --- Ce qui, entre parenthèses, est tellement énorme que ça montre le caractère maladif des réclamations de Julia.
\item[b] Ma démonstration, dit Julia, est entièrement différente de celle de Montel. Montel dit $f(2^nz)$ serait normale or elle ne l'est pas. Je procède en [\string?] inverse et c'est tout différent car Montel utilise les valeurs exceptionnelles pour prouver que la suite n'est pas normale et moi pas. 

Qu'il ait dit qu'il y avait là une différence qu'il tenait à souligner car de là venaient les avantages qu'il tire de sa présentation, rien que de naturel; mais c'est tout autre chose: nous devons juger que ça n'a aucun rapport avec le Montel. Or ça dépend de Montel de deux façons essentielles: les suites normales et puis de la construction de la suite sur laquelle la fonction se reflète, ce qui est la \emph{seule} construction, invention qu'il y ait dans toute [\emph{sic}] ces questions.

Là encore je fais beau jeu à Julia, je ne m'occupe pas de savoir si Montel aurait introduit $f(2^nz)$ si elle n'avait pas été anormale, s'il ne savait pas qu'elle était anormale dans tous les cas, il l'a encore démontré dans la suite du mémoire, s'il ne l'a pas démontré seulement parce qu'il n'en avait pas besoin, si la démonstration d'anormalité est de celles qui comptent ou qui vont de soi, si ce n'est pas de l'ordre de ce qu'on appelle parfois une lacune sans importance, etc. Je donne raison à Julia sur tout cela. Et puis? Sa démonstration reste tributaire des deux faits majeurs: théorie des fonctions normales, construction de la suite $f(2^nz)$.
\item[c] \og Mais j'ai une autre méthode indépendante des familles normales et de la construction de Montel celle qui utilise les résultats de l'Ecole scandinave.\fg Le seul scandinave en l'occurrence est Montel qui démontre un théorème que Lindelöf déclare n'avoir pu ni prouver ni mettre en défaut par un exemple (ce qui n'est tout de même pas une raison pour que ce théorème de Montel ne compte pas à Montel). Et ce théorème n'a pu être prouvé que par les familles normales et les deux méthodes de J. se ressemblent alors comme deux s\oe urs.
\end{enumerate}
/

(4)

Tout long discours demande une conclusion; celui-ci s'en passera ou plutôt, revenant à votre lettre, je conclurai: peut-être que Julia ne dénigre pas l'\OE uvre de Montel mais il ne recule devant rien pour essayer de s'en affranchir. Le résultat ressemble alors si fort à un dénigrement que vous excuserez un myope comme moi s'il s'y est trompé. Pour moi j'aimerais mieux grandir ceux que je voudrais vaincre que de les diminuer.

Voilà, avec les intermédiaires de Rigollot, de ventouses et gargarismes, la journée s'est passée à vous répondre. Mon inoccupation est l'excuse de ma longueur; j'en suis un peu confus mais content de vous avoir dit exactement ma pensée. Je terminerai par un souhait, que votre indécision ne dure plus et que nous puissions dire nettement aux candidats pour qui chacun de nous est.

Sur quoi je m'allonge; je vais certainement bien mieux mais suis faible.

\medskip
à vous

\begin{center}
[signé] H. Lebesgue
\end{center}

\subsection*{La lettre de Montel}

\medskip
\begin{flushright}
Paris, le 24 déc. 1933
\end{flushright}

\medskip
Cher Monsieur Cartan,

\medskip
J'ai trouvé votre lettre hier soir, en rentrant chez moi\footnote{Le 24 décembre 1933 était un dimanche, Montel a trouvé la lettre à laquelle il répond le samedi. Cartan a dû l'écrire le jeudi soir ou le vendredi, après avoir discuté avec Montel le jeudi 21 décembre.}. 

La démonstration de Julia du théorème de Picard est bien identique à celle de la page 299 de mon mémoire de 1916, comme vous me l'avez fait remarquer.

Elle diffère de celle de la page 252 du même mémoire qui s'appuie sur le fait que la fonction est dépourvue de zéro. Mais à cette page 252, je renvoie en Note à une première démonstration du th. de Picard qui se trouve à la page 514 de mon mémoire de 1912. Cette dernière s'applique immédiatement au cas envisagé par Julia.

Pour moi, la question principale / réside dans le fait de substituer au critère des valeurs exceptionnelles un critère quelconque de famille normale et de mentionner que le théorème de Picard s'en déduit.

Je dis que mon théorème de la page 296 va plus loin et comprend le précédent comme cas particulier. Je démontre à cet endroit que tout critère entraîne un nouveau théorème type Landau, ou type Schottky, ou type Picard.

Le raisonnement que je vous ai indiqué jeudi montre que $f(2^nz)$ n peut être normale pour $\module{z}<1$. Pourrait-elle avoir un seul point irrégulier comme dans le cas où $f(z)=z$ par exemple? On voit, comme dans ma démonstration de la page 514 (1912), ou autrement, que cela est en contradiction avec le théorème de Weierstrass. *

Il s'agit dans tout ce qui précède de la démonstration du th. de Picard proprement dit. Il va de soi que je n'ai / jamais songé à revendiquer la notion de droite de Julia. Les indications de mes mémoires ultérieurs sont parfaitement explicites sur ce point et mon opinion est formulée au dernier alinéa de la page 20 de ma Notice.

\medskip
Bien cordialement à vous

\centerline{[signé] Paul Montel}

\medskip
*Mon raisonnement montre en effet qu'aucune suite partielle extraite de $f(2^nz)$ ne peut être normale pour $\module{z}<1$. Donc $f(2^nz)$ tend uniformément vers l'infini pour $\frac{1}{2^2}<\module{z}<\frac{1}{2}$, puisque tout suite partielle tend uniformément vers l'infini.

\section{Commentaires}
\subsection{Pas de lettre de Julia à \'Elie Cartan à ce sujet}
Il n'existe pas à notre connaissance de lettre de  Julia à Cartan à ce sujet et il est assez probable qu'il n'y en a pas eu: Julia habitait à Versailles et Cartan au Chesnay\footnote{Leurs adresses au 1\up{er} février 1934 figurent dans les cahiers \og Vie de la société\fg du \emph{Bulletin} de la \textsc{smf}, 
\begin{itemize}
\item 27 avenue de Montespan, Le Chesnay, pour Cartan,
\item 4 bis rue Traversière, Versailles, pour Julia.
\end{itemize}
Montel, lui, habitait à Paris, 79 rue du Fbg St-Jacques.}, ils ont eu de nombreuses occasions de discuter de vive-voix. Citons, en témoignage, un extrait du discours que fit Julia quelques années plus tard, le 18 mai 1939, à l'occasion du jubilé d'\'Elie Cartan (reproduit dans~\cite{OC6}).

\begin{quote}
Du temps passe; votre ancien élève devenu votre collègue, a dû comme vous émigrer vers la \og cité des eaux\fg. Désormais, c'est presque chaque semaine que le train nous réunit pour nous ramener à Versailles, perdus dans la foule des amateurs de banlieue. Le wagon est presque toujours plein et bruyant, mais l'inconfort ne vous gêne pas, avec votre résistance physique étonnante d'homme issu de la terre. Nous y poursuivons des conversations à bâtons rompus, où les propos universitaires, professionnels ou mathématiques se mêlent aux controverses de tout genre. [...] 

C'est vers cette époque que vous m'avez appelé dans votre maison du Chesnay et dans l'austère petit bureau qu'égayait seulement un carré de fenêtre ouvert sur des branches, et, dans l'ombre, un fauteuil Morris aux bras accueillants. Au cours de fréquents visites, on pouvait là, mieux que dans la salle E\footnote{\`A l'\textsc{ens}, où Julia avait suivi des cours d'agrégation donnés par \'Elie Cartan en 1914, comme il l'évoque au début du même texte.}, mieux que dans les trains de banlieue, mieux que dans les avenues de Versailles, égrener ces propos, coupés de silences, où se révèle l'homme intérieur. [...]
\end{quote}
Julia a donc parlé avec Cartan dans le train, dans la rue à Versailles, et même chez lui au Chesnay où Cartan l'a \og appelé\fg. Le texte permet de dater approximativement cette \og invitation\fg: le musicien Jean Cartan était déjà mort\footnote{Jean Cartan, le deuxième fils d'\'Elie Cartan, un compositeur élève de Dukas et Roussel, est mort de la tuberculose le 26 mars 1932, à l'âge de 25 ans.} (ce fait est évoqué dans les phrases suivantes), d'une part, et d'autre part, Julia n'était pas encore membre de l'Académie des sciences, comme l'indique le paragraphe suivant. Entre 1932 et 1934, pendant la période qui nous intéresse, donc.

\medskip
\'Elie Cartan a laissé la réputation d'un homme d'une grande gentillesse, peut-être un peu indécis, mais juste, \og la conscience même\fg, dit Lebesgue. Il s'est sans doute trouvé là dans une situation fort délicate.

\subsection{Lebesgue et la personnalité de Julia}
Lebesgue reproche à Julia d'être \og âprement personnel\fg. Nous ne savons bien entendu pas ce que Julia a dit à Cartan, mais il est clair qu'il y avait une controverse. Avant de venir à une description de cette controverse, faisons un petit commentaire. Il est très intéressant de noter que Lebesgue envisage, non seulement la blessure de Julia, mais aussi le fait que les aspects qu'il trouve désagréables dans sa personnalité puissent être des conséquences de cette blessure. Par ailleurs, il prend bien garde de dissocier le grand mutilé et ses \og procédés\fg.

En lisant ce qui a été publié sur Julia, notre impression serait plutôt que, bien sûr la blessure inoubliable de Julia est là, nul ne peut l'ignorer, mais on n'en parle pas, sauf peut-être pour le plaindre --- en plus de ses sinistres aspects esthétiques, cette blessure l'a d'ailleurs fait souffrir toute sa vie.

\subsection{Les mathématiques de Montel et de Julia}\label{subsecmathmontel}
La notice que Julia a fait imprimer pour cette élection en 1934 est celle qu'il a incluse dans le premier volume de ses \OE uvres. Julia y est en effet bien moins reconnaissant aux idées de Montel qu'il ne l'avait été en 1917:
\begin{quote}
\`A cette époque, j'ignorais les travaux de M.~Montel. Mon attention sur eux fut attirée par sa Note du 4 juin 1917. Je les étudiai à ce moment dans un tirage à part que M.~Montel voulut bien m'envoyer
\end{quote}
avait-il écrit dans la note~\cite{Julia17} de décembre 1917 (à propos du premier pli cacheté qu'il avait envoyé à l'Académie des sciences sur l'itération). Cette note était une expression de la colère de Julia contre Fatou (voir~\cite{GPSM18}), à un moment où Julia n'était pas (encore) en concurrence avec Montel. \`A propos de ces travaux, Humbert dans son rapport sur le Grand prix des sciences mathématiques avait écrit:
\begin{quote}
\noindent il [Julia] introduit systématiquement, non plus les points invariants attractifs, mais les points invariants où le module du multiplicateur est \emph{supérieur} à l'unité; leur propriété fondamentale est d'être des \emph{points de répulsion}. D'une manière plus précise, si l'on entoure l'un d'eux d'un domaine arbitrairement petit, les conséquents successifs de ce domaine \emph{finissent} par comprendre à leur intérieur tous les points du plan, sauf un ou deux, au plus.

L'analogie de cet énoncé avec celui d'un théorème de M. E. Picard n'a rien de mystérieux: la démonstration de M. Julia repose ici, comme souvent dans le reste du Mémoire, sur la théorie des \emph{suites normales} de M. Montel, théorie dont on sait le lien étroit avec le théorème de M. Picard.
\end{quote}
C'est la note~\cite{Montel17} de Montel qui a déclenché les avancées dans les travaux de Julia et Fatou sur l'itération. Pierre Fatou avait écrit, lui, dans l'introduction du premier de ses articles sur ce sujet~\cite{Fatou19}:
\begin{quote}
\noindent J'ai ensuite reconnu que les recherches récentes relatives aux fonctions analytiques qui admettent des valeurs exceptionnelles permettaient d'aborder le problème dans toute sa généralité; les propriétés des fonctions auxquelles nous faisons allusion sont dues principalement à MM. Landau et Schottky et sont l'extension des théorèmes de M. Picard sur les fonctions entières. M. P. Montel, en leur appliquant la notion de suite normale de fonctions analytiques qui lui est due, est parvenu à des théorèmes élégants et d'un emploi très commode dans les applications. Nous avons fait exclusivement usage dans nos recherches des propositions de M. Montel.
\end{quote}

En 1917, la cause était donc entendue: c'était la notion de famille normale et son utilisation par Montel qui avaient impulsé les travaux sur l'itération. En 1934, dans sa notice, Julia écrit pourtant:
\begin{quote}
En outre, développant logiquement les conséquences de ce principe que les points singuliers d'une fonction la déterminent dans une mesure plus ou moins grande, j'ai recherché la notion correspondante pour les familles et j'ai introduit en Analyse la notion générale de \emph{points singuliers des familles de fonctions d'une ou plusieurs variables}, que d'éminents géomètres\footnote{C'est Ostrowski~\cite{Ostrowski} qui est responsable de la terminologie.} ont bien voulu appeler \emph{points de Julia} ou \emph{points $J$}. Ce concept, aujourd'hui classique, est exposé avec plus ou moins de détails dans les traités de MM.~Goursat, Bieberbach, dans les monographies de MM.~Borel et Valiron sur les fonctions entières, de MM.~Montel et Valiron sur les familles normales.
\end{quote}
Il est certain qu'il s'agit d'un considérable \og amenuisement\fg\footnote{Il est remarquable aussi que toute la notice de Julia soit écrite dans ce style \og personnel\fg... pas au sens où ce style est original mais au sens où Julia revendique, pour son propre compte, il s'attribue, à la première personne du singulier, un certain nombre de qualités et d'idées. \og \^Aprement personnel\fg, disait Lebesgue.}.

\medskip
Le gros \og coup\fg de Montel, c'est son théorème sur les famille normales. Le paragraphe de la lettre de Lebesgue dans lequel celui-ci dit qu'avant Montel \og tout le monde marchait au hasard\fg se réfère à ce théorème. Celui-ci affirme (en termes modernes) que, dans l'espace des fonctions holomorphes  sur un ouvert de $\CC$ (avec la topologie de la convergence uniforme sur les compacts), les compacts sont les fermés bornés.

L'idée nouvelle à laquelle Lebesgue fait allusion lorsqu'il dit:
\begin{quote}
\noindent Des assortiments de théorèmes sont devenus des cas particuliers d'un même fait bien simple.
\end{quote}
c'est d'appliquer un théorème sur les familles, les suites, à l'étude d'une fonction unique, en remplaçant l'étude d'une fonction $f$ en $0$ par l'étude de la suite des fonctions $z\mapsto f(z/2^n)$. Les théorèmes de Picard, par exemple, se démontrent alors comme conséquences du théorème de Montel. 

Rappelons les énoncés de ces théorèmes: 
\begin{itemize}
\item le \og premier théorème de Picard\fg affirme qu'une fonction entière non constante prend toutes les valeurs sauf peut-être une,
\item le \og deuxième théorème de Picard\fg affirme que, si une fonction holomorphe sur le disque épointé a une singularité essentielle en $0$, elle prend toutes les valeurs complexes sauf peut-être une.
\end{itemize}

La démonstration originelle de Picard~\cite{Picard79} (en 1879) utilisait la fonction modulaire. Ce que l'on appelait une \og démonstration élémentaire\fg était une démonstration qui n'en faisait pas usage. La première avait été donnée par Borel~\cite{Borel96} en 1896.

\medskip
Pour en dire un peu plus long sur cette discussion, le plus simple est sans doute de reprendre les choses comme Montel les présente dans sa lettre. Le livre de Julia dont il est question est~\cite{Julia24}. Il s'agit des notes d'un cours Peccot, qu'il a donné en 1920, et qui ont été rédigées par Paul Flamant. Le livre est paru en 1924. Comme il le dit dans sa préface (datée d'août 1923):
\begin{quote}
\noindent [...] l'étude de la fonction dans un cercle entourant le point singulier essentiel peut se ramener à l'étude d'une suite de fonctions dans une couronne circulaire entourant ce point, chacune des fonctions étant \emph{méromorphe} dans la couronne. De ce point de vue,on est conduit, avec M.~Montel, à une nouvelle démonstration des théorèmes de M.~Picard: e Chapitre~III expose les propriétés essentielles des familles normales de fonctions; la fonction modulaire y joue un rôle essentiel pour l'établissement  d'un critérium donné en 1911 par MM.~Landau et Caratheodory sous forme de criterium de convergence, mais dont la transformation en criterium de familles normales permet de redémontrer simplement les théorèmes du Chapitre II.
\end{quote}

C'est des démonstrations de ce livre qu'il est question dans la lettre de Montel.
\begin{itemize}
\item Montel mentionne la démonstration de Julia du théorème de Picard. Il s'agit de celle donnée pages 79-80, dans le chapitre \emph{Les familles normales de fonctions}, qui est en effet identique à la démonstration de Montel, p.~299, sauf que Montel traite le cas plus général d'une fonction méromorphe.
\item Il mentionne aussi sa démonstration, page 252. En effet, il a déjà démontré une première fois el théorème de Picard (pour une fonction analytique, cette fois). Sa démonstration utilise le fait que les $f_n$ de la suite considérée ne prennent pas la valeur $0$. Comme il le dit dans sa lettre, il renvoie en effet à sa démonstration de 1912~\cite[p.~514]{Montel12} (dans laquelle il démontrait le second théorème de Picard, grâce aux familles normales).
\end{itemize}

\section{L'élection de Julia}\label{secelection}
Puisque les élections de Borel et Lebesgue à l'Académie des sciences sont mentionnées dans la lettre de Lebesgue, commençons cette section par une sorte de préhistoire.

\subsection*{L'élection de Borel à l'Académie des sciences}
Les rapports entre Borel et Lebesgue étaient devenus, dès la fin de la guerre, déplorables. La correspondance~\cite{LebesgueBorel,LebesgueBorel2} se termine bien tristement. \`A la rivalité scientifique\footnote{Il est aujourd'hui bien facile de se contenter d'aphorismes, du genre \emph{que serait la mesure de Lebesgue sans les boréliens?}, ou inversement \emph{à quoi bon les boréliens sans théorie de l'intégration?}.} s'ajoute le ressentiment de Lebesgue contre les activités extra-scientifiques auxquelles Borel a consacré de plus en plus de temps au fil des années (la \emph{Revue du mois}, la politique). Les élections que mentionne Lebesgue dans la lettre sont celles 
\begin{itemize}
\item du 11 avril 1921, pour le remplacement de Georges Humbert, le Comité secret (c'est-à-dire les académiciens réunis en séance non publique) avait proposé Borel en première ligne et Lebesgue en deuxième et les 54 votants avaient élu Borel avec 48 voix (et 4 à Lebesgue),
\item du 29 mai 1922, pour le remplacement de Jordan, Lebesgue avait obtenu 44 voix (Vessiot 5, Drach 3, Cartan 2).
\end{itemize}

\subsection*{Avant: l'élection d'\'Elie Cartan}Au moment où les lettres citées ici ont été écrites, le dernier membre de la section de géométrie a avoir été élu est justement \'Elie Cartan. L'élection a eu lieu le 9 mars 1931, il s'agissait alors de remplacer Paul Appell, décédé le 24 octobre 1930. Réunis en comité secret le 2 mars (le lundi précédent l'élection), les géomètres ont classé 
\begin{itemize}
\item en première ligne \'Elie Cartan,
\item en deuxième ligne, ex-aequo, par ordre alphabétique Arnaud Denjoy, Gaston Julia, Paul Montel et Ernest Vessiot.
\end{itemize}
Les cinquante-quatre académiciens présents lors de l'élection ont ensuite suivi l'avis de la section, Cartan obtenant 51 voix, Vessiot 2 et Montel 1.

\medskip
Nous l'avons dit, Painlevé est mort le 29 octobre 1933. L'Académie des sciences va donc procéder à son remplacement. Toutes les informations qui suivent viennent, soit des \emph{Comptes rendus}, soit du registre des Comités secrets (aux archives de l'Académie des sciences).

\subsection*{Le remplacement de Painlevé}La procédure habituelle se déroule pour le remplacement de Painlevé. Le 19 février 1934, Julia a annoncé officiellement sa candidature, le 26 février, la section, constituée de Jacques Hadamard, \'Edouard Goursat, \'Emile Borel, Henri Lebesgue et \'Elie Cartan\footnote{Picard, depuis qu'il a été élu Secrétaire perpétuel en 1917, n'est plus membre de la section (il a été remplacé par Goursat). Il participe néanmoins, en tant que Secrétaire perpétuel, à la réunion.}, se réunit en comité secret. La réunion dure près de deux heures\footnote{Les heures du début et de la fin de la réunion sont indiquées dans les \emph{Comptes rendus}.}, ce qui est exceptionnellement long et un signe que les choix n'étaient pas simples. Le registre des comités secrets contient des indications des positions des uns et des autres:
\begin{itemize}
\item \'Elie Cartan expose les raisons pour lesquelles il votera pour Montel --- il fait également l'éloge des travaux de Julia,
\item Picard place Julia en première ligne,
\item Lebesgue place Montel en première ligne. 
\end{itemize}
Hadamard, Goursat et Borel interviennent eux aussi dans la discussion. 

\begin{remarque*}
Il est vraisemblable qu'Hadamard a lui aussi voté pour Montel, comme l'indique une discussion qui a eu lieu quelque temps auparavant. Le 14 novembre~1932, l'Académie des sciences a discuté en comité secret du prix Albert I\up{er} de Monaco et en particulier, à la suite de Picard, s'est demandée à quel géomètre elle pourrait décerner ce prix. Il a été question de Julia, de Montel, de Vessiot et de Denjoy. Lebesgue a  déclaré que, s'il s'agissait d'une élection dans la section de géométrie, il voterait pour Montel et Hadamard a exprimé la même opinion\footnote{Le prix Albert I\up{er} a finalement été décerné à Louis de Broglie.}.
\end{remarque*}

Finalement, la section a proposé
\begin{itemize}
\item en première ligne, Paul Montel,
\item en deuxième ligne, Gaston Julia,
\item en troisième ligne, Arnaud Denjoy,
\item en quatrième ligne, ex-aequo, Maurice Fréchet, René Garnier et Paul Lévy.
\end{itemize}

Le lundi suivant, qui était le 5 mars 1934, les académiciens des sciences se sont prononcés. Les cinquante-quatre présents ont élu Julia, avec 32 voix contre 21 à Montel (et un bulletin nul). En d'autres termes, ils n'ont pas suivi l'avis de la section\footnote{Ajoutons que le délai entre la mort de Painlevé et son remplacement, du 29 octobre au 5 mars, un peu plus de quatre mois, est exactement le même que celui écoué entre la mort de Jordan et l'élection de Lebesque (voir la note~\ref{notemortjordan}), qui n'avait donc rien d'exceptionnellement long.}.

\medskip
Un tel renversement (en tout cas lors de l'élection d'un géomètre) est assez rare\footnote{Pendant la période de l'entre-deux-guerres, l'Académie des sciences a recruté plusieurs membres pour la section de mécanique, de façon un peu plus mouvementée:
\begin{itemize}
\item Jules Drach en 1929, qui avait été rajouté (sur proposition de Lebesgue, d'après le registre des comités secrets, 3 juin 1929) par l'Académie à la liste proposée par la section et élu avec une majorité confortable,
\item Jouguet, Villat et Louis de Broglie en 1930, 1932, 1933, sans problème,
\item Alfred Caquot en 1934, passé devant Vessiot (renversement de la proposition de la section) et au second tour.
\end{itemize}}. Dans les cinquante années précédentes, il y a eu treize élections de géomètres et une seule fois, les académiciens n'ont pas suivi la section, le 19~mai 1919, où ils ont préféré Goursat à Borel, avec un score assez semblable: 29 pour Goursat soutenu, d'après le registre des comités secrets, par Picard et Jordan, contre 23 pour Borel, soutenu par Humbert et Painlevé. Dans la plupart des cas, comme dans celui de l'élection de Cartan, la majorité obtenue par l'heureux élu est très confortable (Darboux obtient 47 voix sur 53, Appell 52 sur 53, Humbert 54 sur 58...) mais il y a eu des élections plus difficiles, comme celle de... Poincaré (qui a obtenu 31 voix contre 24 à Mannheim, un géomètre un peu oublié aujourd'hui) le 31 janvier 1887, ou même celle d'Hadamard (36 voix contre 21 à Goursat)\footnote{Picard soutenait déjà Goursat contre Painlevé en 1900 (voir la lettre de Painlevé à Mittag-Leffler du 10 juin 1900 citée dans~\cite[p.~811]{PainleveOC3}). Après Painlevé, Humbert et Hadamard ont été préférés à Goursat par les académiciens. Les clivages politiques ont joué leur rôle ici aussi. Dans le cas de Jacques Hadamard, il est assez vraisemblable que, quelques années après la conclusion de l'Affaire Dreyfus, les anti-dreyfusards ont voté contre un candidat qui, non content d'être juif, avait été un dreyfusard très actif et était toujours militant de la \og Ligue des droits de l'Homme\fg.}.

\begin{remarque*}
Il est certain que le fait que Julia ait été un grand blessé de la grande guerre a grandement joué en sa faveur. Voici quelques exemples de situations analogues:
\begin{itemize}
\item le 26 mai 1919, lors d'une élection pour la section de chimie, Hadamard explique qu'il \og votera pour M. Béhal à cause du rôle de premier plan qu'il a joué dans les fabrications chimiques de guerre\fg,
\item le 4 avril 1921, les services rendus par \'Emile Borel à la défense nationale pendant la guerre ont été exposés en faveur de son élection,
\item trois notes publiées dans les \emph{Comptes rendus} en 1916 par Esclangon\footnote{Sous les titres 
\begin{itemize}
\item Sur les trajectoires aériennes des projectiles, 
\item Sur les coups de canon et les zones de silence,
\item Sur le principe de Doppler et le sifflement des projectiles.
\end{itemize}} seront utiles lors de son élection en 1929~\cite[p.~139]{saintmartin08}, 
\item encore en 1942, lorsqu'il s'agira de remplacer Lebesgue\index{Lebesgue (Henri), 1875--1941, mathématicien}, un académicien écrira à Fréchet\index{Frechet@Fréchet (Maurice), 1878--1973, mathématicien} que le fait qu'il soit ancien combattant jouera en sa faveur (lettre citée dans~\cite[p.~365]{TaylorFrechet2}).
\end{itemize}
\end{remarque*}

\subsection*{Remplacement de Goursat}
L'élection suivante pour la section de géométrie a eu lieu le 31 mai 1937, pour remplacer Goursat décédé le 25 novembre 1936. La réunion en comité secret du 24 mai n'a duré que quinze minutes et a établi la liste
\begin{itemize}
\item en première ligne, Paul Montel,
\item en deuxième ligne, Arnaud Denjoy,
\item en troisième ligne, ex-aequo, Maurice Fréchet, René Garnier, Paul Lévy et Georges Valiron,
\end{itemize}
et le vote donnera 51 voix à Montel, 1 à Garnier, 1 à Paul Lévy... et rien à Denjoy, qui sera le prochain élu, en 1942, bien avant Garnier (1952), Fréchet (1956) et Lévy (1964).

La cause était entendue. Denjoy, dont Lebesgue nous a dit dans sa lettre qu'il s'était retiré devant Montel lors de l'élection de 1934, n'avait d'ailleurs même pas actualisé sa notice pour cette élection en 1937 parce qu'\og il y avait unanimité des académiciens pour Montel\fg (dossier Denjoy, archives de l'Académie des sciences).

\section*{Remerciements}
Nous remercions la famille Cartan pour l'autorisation de publier ces deux lettres.

\nocite{mathnicois}

\newcommand{\SortNoop}[1]{}
\providecommand{\bysame}{\leavevmode ---\ }
\providecommand{\og}{``}
\providecommand{\fg}{''}
\providecommand{\smfandname}{\&}
\providecommand{\smfedsname}{\'eds.}
\providecommand{\smfedname}{\'ed.}
\providecommand{\smfmastersthesisname}{M\'emoire}
\providecommand{\smfphdthesisname}{Th\`ese}

\vfill

\end{document}